%% file: enu_vol_surv.tex
\documentclass [11pt,oneside]{article}

\usepackage{amssymb}
\usepackage{amsfonts}
\usepackage{amsmath}
\usepackage{amsthm}
\usepackage{epsfig}

\newtheorem{lemma}{Lemma}[section]
\newtheorem{teo}[lemma]{Theorem}
\newtheorem{rem}[lemma]{Remark}
\newtheorem{prop}[lemma]{Proposition}
\newtheorem{cor}[lemma]{Corollary}

\newcommand{\matE} {\ensuremath {\mathbb{E}}}
\newcommand{\matM} {\ensuremath {\mathbb{M}}}
\newcommand{\matR} {\ensuremath {\mathbb{R}}}

\newcommand{\matZ} {\ensuremath {\mathbb{Z}}}
\newcommand{\matC} {\ensuremath {\mathbb{C}}}

\newcommand{\matH} {\ensuremath {\mathbb{H}}}
\newcommand{\matB} {\ensuremath {\mathbb{B}}}

\newcommand{\calK} {\ensuremath {\mathcal{K}}}

\newcommand{\calM} {\ensuremath {\mathcal{M}}}

\newcommand{\calA} {\ensuremath {\mathcal{A}}}

\newcommand{\calD} {\ensuremath {\mathcal{D}}}

\newcommand{\calP} {\ensuremath {\mathcal{P}}}
\newcommand{\calT} {\ensuremath {\mathcal{T}}}

\newcommand{\calH} {\ensuremath {\mathcal{H}}}
\newcommand{\calHtilde} {\ensuremath {\widetilde{\mathcal{H}}}}

\newcommand{\TV}{{\rm TV}}

\newcommand{\mylegend} [1] {\caption{\footnotesize{#1}}}

\newfont{\Got}{eufm10 scaled 1200}

\newcommand{\mynote}[1]{}
\font\titsc=cmcsc10 scaled 1200

\newcommand{\Ybar}{\overline{Y}}
\newcommand{\Ybarbar}{\overline{\overline{Y}}}
\newcommand{\Ptilde}{\widetilde{P}}

\newcommand{\mettifig}[1]{\epsfig{file=#1}}

\author{Alexander \titsc{Mednykh} \and Carlo \titsc{Petronio}\thanks{Both
authors were supported by the INTAS Project ``CalcoMet-GT'' 03-51-3663}}

\title{Hyperbolic 3-manifolds with geodesic boundary:\\ Enumeration and volume calculation}

\begin{document}

\maketitle

\begin{abstract}
    \noindent
    We describe a natural strategy to enumerate compact hyperbolic $3$-manifolds with
    geodesic boundary in increasing order of complexity. We show that the same strategy
    can be employed to analyze simultaneously compact manifolds and finite-volume
    manifolds having toric cusps.  In opposition to this we show that,
    if one allows annular cusps, the number of manifolds grows very rapidly,
    and that our strategy cannot be employed to obtain a complete list.
    We also carefully describe how to compute the volume of our manifolds,
    discussing formulae for the volume of a tetrahedron with generic dihedral angles in
    hyperbolic space.
  \vspace{4pt}

\end{abstract}

\noindent According to Thurston's geometrization program, the
theory of hyperbolic manifolds plays a central role in
$3$-dimensional topology. Hyperbolic manifolds with geodesic
boundary, the first example of which was given by Thurston himself
in~\cite{Thurston:notes} (and later generalized in~\cite{PaZi}),
are an important portion of this theory. On the other hand, the
algorithmic and computer approach to $3$-manifolds has been
acquiring an increasing popularity in recent years. For cusped
hyperbolic manifolds this approach, which was worked out
in~\cite{CaHiWe, SaWe, Weeks:tilt, snappea} and several other
papers, again depends on ideas of Thurston, namely on the use of
moduli for hyperbolic ideal tetrahedra and equations to ensure
consistency of the structures. For closed manifolds the basics of
the computer approach were set by Matveev
in~\cite{Matveev:complexity} (see also~\cite{FoMa,
Matveev:ATC3M}), and several experimental results were later
obtained by himself and other authors 
(see~\cite{Bruno:surv, MaPe}). In the present paper we describe the general setting
of the algorithmic approach to hyperbolic $3$-manifolds with
geodesic boundary, concentrating in particular on their
enumeration in order of increasing complexity, and on the
computation of their volume.

\section{Hyperbolic structures}\label{structure:section}
In this section we review the general theory of hyperbolic
$3$-manifolds with geodesic boundary, stating the main results we
will need in the sequel.

\paragraph{Local structure, cusps, compactifications}
In the rest of this paper we will call \emph{hyperbolic} an
orientable finite-volume complete Riemannian $3$-manifold with
\emph{non-empty} boundary, locally isometric to an open subset of
a closed half-space of $\matH^3$. We will always denote such a
manifold by $Y$. Note that $\partial Y$ is totally geodesic.
Doubling $Y$ along its boundary and using the description of the
ends of the double~\cite{BePe:book}, one can show that $Y$
consists of a compact portion together with some ``toric and
annular cusps.'' A toric cusp is here a space of the form
$T\times[0,\infty)$, attached to the rest of $Y$ along $T\times
\{0\}$, and an annular cusp is defined analogously. Note that
toric cusps are disjoint from $\partial Y$, while an annular cusp
gives two punctures in $\partial Y$. In particular, $\partial Y$
is compact if and only if $Y$ has no annular cusps.

Given $Y$ as above, we can naturally get a compact manifold
$\Ybar$ by adding a torus $T\times\{\infty\}$ for each toric cusp
$T\times[0,\infty)$, and an annulus $A\times\{\infty\}$ for each
annular cusp $A\times[0,\infty)$. However, it turns out that when
there are annular cusps another compactification $\Ybarbar$ of $Y$
is more suited to the geometric situation. We define $\Ybarbar$ as
a quotient of $\Ybar$, where in each annular cusp
$A\times[0,\infty]$ with $A=S^1\times[0,1]$, for all $t\in[0,1]$
we collapse $S^1\times \{t\}\times\{\infty\}$ to one point
$\{*\}\times\{t\}\times\{\infty\}$. Note that $Y$ is obtained from
$\Ybarbar$ by removing some toric boundary components and drilling
some properly embedded arcs.

\paragraph{Rigidity and Kojima decomposition}
A key result for computational purposes is the following:

\begin{teo}\label{rigidity:teo}
Any two homeomorphic hyperbolic manifolds are isometric.
\end{teo}

This result is commonly referred to as \emph{rigidity theorem},
and a proof for the geodesic boundary case was spelled out
in~\cite{FriPe}.

Another very important fact is that the hyperbolic structure
determines certain combinatorial data which can be employed to
efficiently test two manifolds for homeomorphism. This result is
analogous to the Epstein-Penner decomposition of cusped hyperbolic
manifolds without boundary~\cite{EpPe}, and it was proved by
Kojima~\cite{Kojima:poly, Kojima:poly-bis}. Its statement involves
the notion of \emph{truncated polyhedron}, that we now give.
Consider the projective model of hyperbolic $3$-space, viewed as
the open unit ball $\matB^3$ in the Euclidean $3$-space $\matE^3$.
Let us call \emph{finite} the points of $\matB^3$, \emph{ideal}
those of $\partial\matB^3$, and \emph{ultra-ideal} the other
points of $\matE^3$. Consider a convex polyhedron $\Ptilde$ in
$\matE^3$ with ideal and/or ultra-ideal vertices, and all edges
meeting the closure of $\matB^3$. Dual to each ultra-ideal vertex
of $\Ptilde$ there is an open hyperbolic half-space, and we define
$P$ to be $\Ptilde\cap\matB^3$ minus these half-spaces. Any $P$
arising like this from some $\Ptilde$ will be called a truncated
polyhedron. Note that $P$ has \emph{internal} faces, those coming
from faces of $\Ptilde$, and \emph{truncation} faces. Moreover
internal and truncation faces lie at right angles to each other.

\begin{teo}\label{kojima:teo}
Any hyperbolic manifold admits a canonical decomposition as a gluing
of truncated polyhedra along the internal faces.
\end{teo}

In the sequel we will need to refer to the geometric argument
underlying this result, so we briefly sketch it here. Regard
$\matE^3$ as the hyperplane at height $1$ in Minkowsky $4$-space
$\matM^{\,3,1}$. Suppose first that $Y$ has no toric cusps, and
identify the universal cover of $Y$ to an intersection of
half-spaces of $\matB^3$.  Dual to each such subspace there is a
point having norm $1$ in $\matM^{\,3,1}$, and the decomposition of
$Y$ is obtained by taking the faces of the convex hull of all
these points, projecting first to $\matE^3$, truncating, and then
projecting to $Y$. When there are toric cusps one must also take
suitably small Margulis neighbourhoods of these cusps, lift them
to horoballs in $\matB^3$, consider the duals to these horoballs
on the light-cone of $\matM^{\,3,1}$, include these duals in the
convex hull, and suitably subdivide some of the resulting faces.
We address the reader to~\cite{FriPe} for all the details.

We are now in a position to explain why we have introduced the
compactification $\Ybarbar$:

\begin{prop}\label{barbar:prop}
Let $Y$ be hyperbolic. For each truncated polyhedron $P$ in the Kojima
decomposition of $Y$, consider the corresponding Euclidean
polyhedron $\Ptilde$. Glue these $\Ptilde$'s along the same maps as in the
Kojima decomposition, and remove open stars of the vertices.  The resulting
space is then  homeomorphic to $\Ybarbar$.
\end{prop}

\paragraph{Topological obstructions to hyperbolicity}
Since dealing with open manifolds is impossible by computer, the
idea to enumerate the hyperbolic $Y$'s is to enumerate the
corresponding compact $\Ybar$'s or $\Ybarbar$'s. The next two
results show that the presence of annular cusps makes a dramatic
difference. The former can be found in~\cite{FriPe}, the latter
follows from the results discussed in Section~\ref{CFMP:section}.
(Recall that we are calling `hyperbolic' a manifold with
\emph{non-empty} geodesic boundary).

\begin{teo}\label{obstructions:teo}
Let $M$ be a compact orientable $3$-manifold with boundary, and let $Y$ be $M$
minus the toric components of $\partial M$. The following conditions are
pairwise equivalent:
\begin{itemize}
\item $Y$ is hyperbolic; \item $Y$ is hyperbolic, it has no
annular cusps, and $\Ybar\cong M$; \item $M$ is irreducible,
boundary-irreducible, acylindrical and atoroidal, and
$\chi(M)<0$.
\end{itemize}
\end{teo}

\begin{prop}\label{bad:annuli:prop}
\begin{itemize}
\item For any $g\geqslant 2$ there exists a hyperbolic $Y$ such that $\Ybar$ is the
handlebody of genus $g$;
\item For any compact $3$-manifold $M$ there exists a hyperbolic $Y$ such that $\Ybarbar\cong M$.
\end{itemize}
\end{prop}

\section{Complexity of 3-manifolds\\
and enumeration strategy}\label{complexity:section}

In this section we recall the basics of the theory of simple and
special spines, and the related theory of ideal triangulations,
describing how it can be employed to (partially) enumerate the
class of $3$-manifolds we are interested in. Proofs of all results
on spines and complexity can be found in~\cite{Matveev:ATC3M}.

\begin{figure}
    \begin{center}
    \mettifig{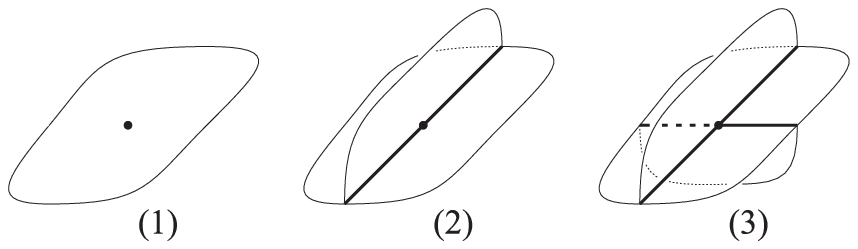,width=9cm}
    \mylegend{Local aspect of an almost-special polyhedron.} \label{almostspecial:fig}
    \end{center}
    \end{figure}

\paragraph{Simple spines and complexity}
Throughout the present section we will employ the PL category for
$3$-manifolds and use the customary notions of PL topology,
see~\cite{RoSa}. A \emph{simple polyhedron} is a compact
polyhedron $P$ such that the link of each point of $P$ can be
embedded in the space given by a circle with three radii. In
particular, $P$ has dimension at most $2$. Finite graphs and
closed surfaces are examples of simple polyhedra. A point of a
simple polyhedron is called a \emph{vertex} if its link is
precisely given by a circle with three radii. A regular
neighbourhood of a vertex is shown in
Fig.~\ref{almostspecial:fig}-(3).
From the figure one sees that the vertices are isolated, whence
finite in number.  Graphs and surfaces do not contain vertices.

If $M$ is a compact $3$-manifold with non-empty boundary, we call
\emph{spine} of $M$ a subpolyhedron $P$ of $M$ such that
$M\setminus P$ is an open collar of $\partial M$. We call
\emph{complexity} of $M$, and denote by $c(M)$, the minimal number
of vertices of a simple spine of $M$. We say that a spine of $M$
is \emph{minimal} if it has $c(M)$ vertices and it does not
contain any proper subpolyhedron which is also a spine of $M$.

\paragraph{Special spines and ideal triangulations}
We now introduce a more restrictive type of spine which turns out
to have a very clear geometric counterpart. A simple polyhedron
$P$ is called \emph{almost-special} if the link of each point of
$P$ is given by a circle with either zero, or two, or three radii.
The local aspects of $P$ are correspondingly shown in
Fig.~\ref{almostspecial:fig}. The points of type (2) or (3) are
called \emph{singular}, and the set of singular points of $P$ is
denoted by $S(P)$. We will say that $P$ is \emph{special} if it is
almost-special, $S(P)$ contains no circle component, and
$P\setminus S(P)$ consists of open $2$-discs.

\begin{figure}
    \begin{center}
    \mettifig{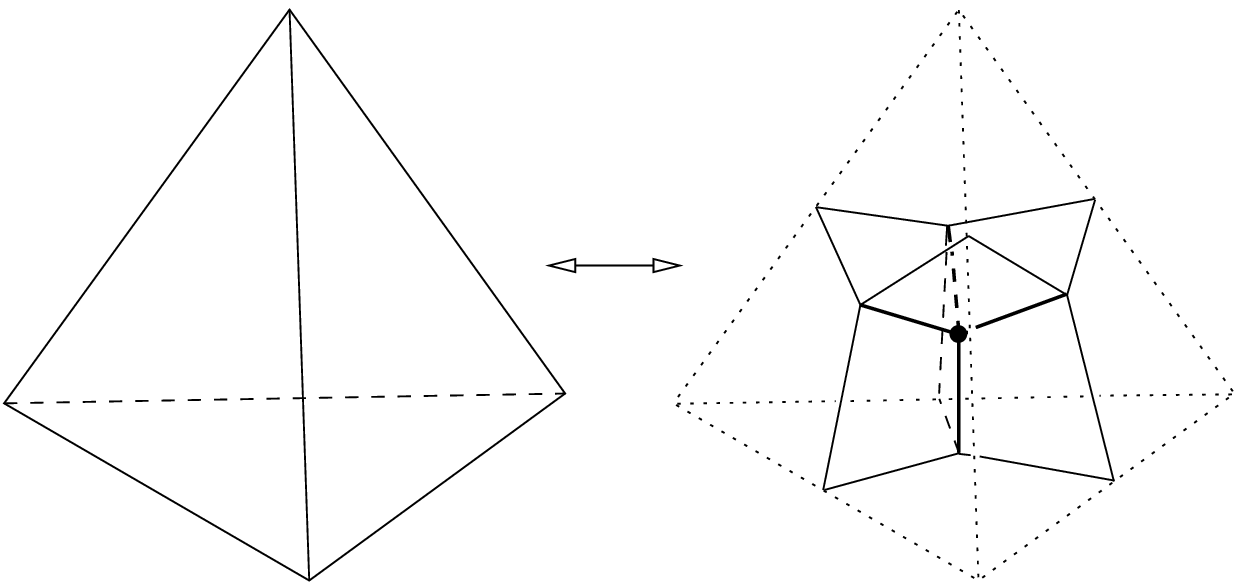,width=8cm}
    \mylegend{Duality between triangulations and special polyhedra.} \label{duality:fig}
    \end{center}
    \end{figure}

We now call \emph{ideal triangulation} of a compact $3$-manifold
$M$ with non-empty boundary a realization of the interior of $M$
as follows. We take a finite number of tetrahedra, we glue
together in pairs the faces of these tetrahedra along simplicial
maps, and we remove the vertices. Equivalently, an ideal
triangulation of $M$ is a realization of $M$ as a gluing of
\emph{truncated} tetrahedra. The relation between spines and
triangulations is given by the following:

\begin{prop}\label{duality:basic:prop}
The set of ideal triangulations of a $3$-manifold $M$ corresponds
bijectively to the set of special spines of $M$. The polyhedron
corresponding to a triangulation is the $2$-skeleton of the dual
cellularization, as shown in Fig.~\ref{duality:fig}.
\end{prop}

\paragraph{Manifolds having special minimal spines}
Special spines have two main advantages if compared to merely
simple ones. First of all, a special spine determines the manifold
it is a spine of~\cite{Casler}, which is false for simple spines.
Second, no efficient method for listing simple spines is known,
whereas enumerating special spines in increasing order of
complexity is very easy (at least theoretically,
see~\cite{BePe:manus}). For these reasons, the next result is
crucial for computational purposes:

\begin{prop}\label{special:prop}
Let $M$ be a compact $3$-manifold with non-empty boundary. The
following conditions are pairwise equivalent:
\begin{itemize}
\item $M$ is irreducible, $\partial$-irreducible, acylindrical,
and not the $3$-disc;
\item $M$ has some special minimal spine;
\item All minimal spines of $M$ are special.
\end{itemize}
\end{prop}

This result is not quite stated in this form
in~\cite{Matveev:ATC3M} or in any of Matveev's papers, but it
easily follows from the proof of~\cite[Theorem
2.2.4]{Matveev:ATC3M}.

\paragraph{Na\"\i f enumeration strategy}
Cusped hyperbolic manifolds without boun\-dary were studied and
enumerated in~\cite{CaHiWe}, which explains why we have decided to
restrict to manifolds with non-empty boundary. In addition, we
forbid here the presence of annular cusps, because, according to
Propositions~\ref{bad:annuli:prop} and~\ref{special:prop}, the
theory of spines, triangulations, and complexity does not appear
to be well-suited for the investigation of such manifolds. See
Section~\ref{CFMP:section}.

Let us then define $\calH_n$ as the set of all hyperbolic
manifolds having complexity $n$ and non-empty \emph{compact}
geodesic boundary. We also define $\calHtilde_n$ as the set of all
orientable, compact, irreducible, $\partial$-irreducible, and
acylindrical manifolds with negative $\chi$.
Theorem~\ref{obstructions:teo} implies that if $Y\in\calH_n$ then
$\Ybar\in\calHtilde_n$. Conversely, if $M\in\calHtilde_n$, then $M$
minus the toric components of $\partial M$ belongs to $\calH_n$ if
and only if $M$ is atoroidal. We can then view $\calHtilde_n$ as
the set of \emph{candidate hyperbolic} manifolds with complexity
$n$ and without annular cusps.

Note now that Proposition~\ref{special:prop} applies to the
elements of $\calHtilde_n$. Therefore the following theoretical
steps lead to the exact determination of $\calH_n$, assuming
$\calH_{m}$ is known inductively for $m<n$:

\begin{enumerate}

\item Produce the list of all special spines with $n$ vertices and
negative $\chi$;

\item Remove from the list the spines $P$ whose associated
manifold $M(P)$ is not hyperbolic;

\item Remove from the list the spines $P$ such that $M(P)$ belongs
to some $\calH_m$ for $m<n$;

\item For $P$ varying in the list, compare the manifolds $M(P)$
for equality, discarding duplicates.

\end{enumerate}

Step (1) does not present any theoretical difficulty, but its
practical implementation is quite demanding if no computational
shortcuts are employed. We will discuss these shortcuts in the
next paragraph. And we will explain how to carry out the other
steps in the next section.

\paragraph{Pseudo-minimal spines}
Taking into account Propositions~\ref{duality:basic:prop}
and~\ref{special:prop}, the reader may wonder why we have employed
special spines rather than triangulations to list the elements of
$\calH_n$. The first remark is that one can always restrict to
\emph{minimal} triangulations without losing any potentially
interesting manifold. However, it turns out that dual to a
triangulation which is minimal among triangulations of the same
manifold, there is often a special spine which is not minimal
among simple spines of the same manifold. Of course this can only
happen if the corresponding manifold violates one of the
topological constraints of Proposition~\ref{special:prop}, but we
know from Theorem~\ref{obstructions:teo} that in this case the
manifold is not hyperbolic, so we can discard it. In other words,
\emph{non-minimality is a much more flexible notion for spines
than for triangulations}, so using spines we can substantially
reduce the list of manifolds that we will later need to
investigate.

\begin{figure}
    \begin{center}
    \mettifig{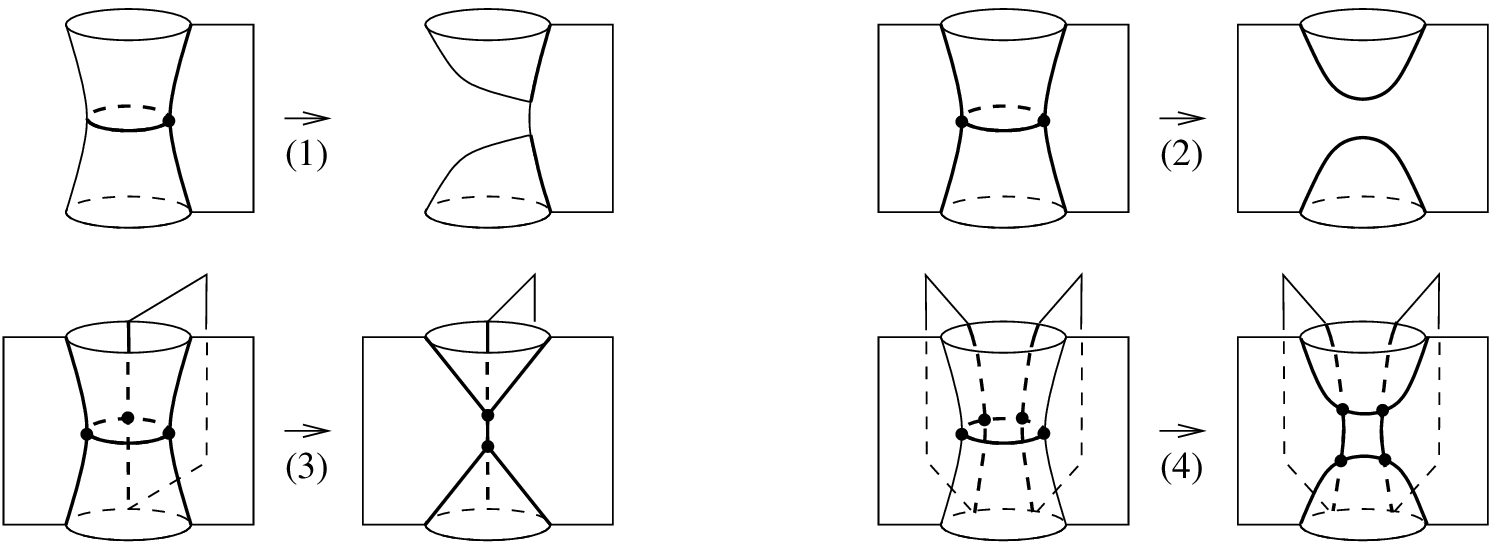, width=12 cm}
    \mylegend{Moves on simple spines.} \label{moves:fig}
    \end{center}
    \end{figure}

Of course it is impossible to check in a direct fashion whether a
special spine is minimal, but there are many criteria for
non-minimality, which can be used as tests to discard spines which
will certainly not bring any relevant manifold. The tests used
in~\cite{FriMaPe3} are based on the moves shown in
Fig.~\ref{moves:fig}, which are easily seen to transform a spine
of a manifold into another spine of the same manifold. To be
precise, let us say that a spine is \emph{pseudo-minimal} if it
cannot be transformed into a spine with fewer vertices by a
combination of the moves shown in Fig.~\ref{moves:fig}. The key
point mentioned above on the flexibility of spines is that moves
(1) and (2) do not lead to special spines (in general), so they do
not have counterparts at the level of triangulations. The first
step of the enumeration strategy is then replaced by the
following:

\begin{itemize}

\item Produce the list of all pseudo-minimal special spines with
$n$ vertices and negative $\chi$.

\end{itemize}

We also mention that another trick very important for
computational purposes is to construct the candidate spines
portion after portion, following the branches of a tree, and to
apply the non-minimality tests arising from the moves of
Fig.~\ref{moves:fig} also to partially constructed spines, thus
``cutting the dead branches'' at their bases.

\section{Hyperbolicity equations and tilts}\label{equations:section}
In this section we discuss how an ideal triangulation can be
employed to construct a hyperbolic structure on a given manifold
and to recognize the canonical Kojima decomposition of that
manifold. This allows to carry out steps (2)-(4) of the
enumeration strategy for $\calH_n$ explained in the previous
section. We actually include in the discussion also manifolds with
annular cusps, because the methods via ideal triangulations to
construct and recognize the hyperbolic structure apply to these
manifolds too. It is only step (1) of the enumeration strategy
(the listing of special spines) that is not suited to manifolds
with annular cusps, and the reason why these manifold are ruled
out from $\calH_n$.

We first treat the compact case and then sketch the variations
needed for the case where there are also some cusps. For all
details and proofs (and for some very natural terminology that we
use here without giving actual definitions) we address the reader
to~\cite{FriPe}.

\paragraph{Compact case: moduli and equations}
The basic idea for constructing a hyperbolic structure via an
ideal triangulation is to realize the tetrahedra of the
triangulation as truncated tetrahedra in $\matH^3$ and then
require that the structures match when the tetrahedra are glued
together. The following facts show that one can use the dihedral
angles as moduli to parameterize the realizations of a tetrahedron
and to check consistency:

\begin{itemize}

\item A hyperbolic structure on a combinatorial truncated
tetrahedron is determined by the 6-tuple of dihedral angles along
the internal edges;

\item The only restriction on this 6-tuple of positive reals is
that the angles of each of the four truncation triangles should
sum up to less than $\pi$;

\item The lengths of the internal edges can be computed as
explicit functions of the dihedral angles;

\item A choice of hyperbolic structures on the tetrahedra of an
ideal triangulation of a manifold $M$ gives rise to a hyperbolic
structure on $M$ if and only if all matching edges have the same
length and the total dihedral angle around each edge of $M$ is
$2\pi$.

\end{itemize}

Given an ideal triangulation consisting of $n$ tetrahedra one then
has the \emph{hyperbolicity equations}: a system of $6n$ equations
with unknown varying in an open set of $\matR^{6n}$ which, by
rigidity (Theorem~\ref{rigidity:teo}), admits one solution at
most.

\paragraph{Canonical decomposition and tilts}
Once a hyperbolic structure has been constructed on a manifold $Y$
using an ideal triangulation $\calT$, one natural issue is to
decide if $\calT$ is the canonical decomposition of $Y$ and, if
not, to promote $\calT$ to become canonical. These matters are
faced using the \emph{tilt formula}~\cite{Ushijima:unified,
Ushijima:tilt, Weeks:tilt}, that we now describe.

Recall first that the Kojima decomposition of $Y$ is constructed
by projecting first to $\matH^3$ and then to $Y$ the faces of the
convex hull of the set $\calP\subset \matM^{\,3,1}$ of the duals
to the boundary components of the universal cover of $Y$. If
$\sigma$ is a $d$-simplex in $\calT$, each end of a lifting of
$\sigma$ to $\matH^3$ determines a point of $\calP$. Now let two
tetrahedra $\Delta_1$ and $\Delta_2$ share a $2$-face $F$, and let
$\widetilde{\Delta}_1,\widetilde{\Delta}_2$ and $\widetilde F$ be
liftings of $\Delta_1,\Delta_2$ and $F$ to $\matH^3$ such that
$\widetilde{\Delta}_1\cap\widetilde{\Delta}_2=\widetilde{F}$. Let
$\overline{F}$ be the $2$-subspace of $\matM^{\,3,1}$ that
contains the three points of $\calP$ determined by
$\widetilde{F}$. For $i=1,2$ let $\overline{\Delta}^{(F)}_i$ be
the half-$3$-subspace bounded by $\overline{F}$ and containing the
fourth point of $\calP$ determined by $\widetilde{\Delta}_i$. Then
one can show that $\calT$ is the canonical Kojima decomposition of
$Y$ if and only if, whatever $F,\Delta_1,\Delta_2$, the following
conditions are met:

\begin{itemize}

\item[(a)] the convex hull of $\overline{\Delta}^{(F)}_1$ and
$\overline{\Delta}^{(F)}_2$ does not contain the origin of
$\matM^{\,3,1}$;

\item[(b)] $\overline{\Delta}^{(F)}_1$ and
$\overline{\Delta}^{(F)}_2$ lie on distinct $3$-subspaces of
$\matM^{\,3,1}$.

\end{itemize}
In addition, if condition (a) is met for all triples
$F,\Delta_1,\Delta_2$, the canonical decomposition is obtained by
merging together the tetrahedra along which condition (b) is not
met.

The tilt formula defines a real number $t(\Delta,F)$ describing
the ``slope'' of $\overline{\Delta}^{(F)}$. More precisely, one
can translate conditions (a) and (b) into the inequalities
$t(\Delta_1,F)+t(\Delta_2,F)\leqslant 0$ and
$t(\Delta_1,F)+t(\Delta_2,F)\neq 0$ respectively. Since we can
compute tilts explicitly in terms of dihedral angles, this gives a
very efficient criterion to determine whether $\calT$ is canonical
or a subdivision of the canonical decomposition. Even more, it
suggests where to change $\calT$ in order to make it more likely
to be canonical, namely along $2$-faces where the total tilt is
positive.  This is achieved by 2-to-3 moves along the offending
faces, as discussed in~\cite{FriPe}.

\paragraph{The non-compact case}
When one is willing to accept both compact geodesic boundary and
cusps, the same strategy for constructing the structure and
finding the canonical decomposition applies, but many subtleties
and variations have to be taken into account. Let us quickly
mention which.

\emph{Moduli}. As suggested by Proposition~\ref{barbar:prop}, to
construct a hyperbolic structure on a manifold $Y$ one must take
an ideal triangulation $\calT$ of the compactification $\Ybarbar$,
in which each toric cusp is completed with a torus and each
annular cusp is completed with a segment. Moreover one must
suppose that the segments arising from the annular cusps of $Y$
are edges of $\calT$. To parameterize tetrahedra one should then
proceed as follows:

\begin{itemize}

\item Assign dihedral angle $0$ to each edge corresponding to an
annular cusp of $Y$, which geometrically means that the edge,
before truncation, is tangent to the boundary of $\matB^3$;

\item If three edges meet at a vertex asymptotic to a toric cusp
of $Y$, assign them dihedral angles summing up to $\pi$, which
geometrically means that the vertex is an ideal one.

\end{itemize}

Switching viewpoint, if one starts from an ideal triangulation $T$
of a manifold $M$, one must arbitrarily choose a family $\alpha$
of edges of $T$ and assign dihedral angles $0$ to the edges in
$\alpha$, and dihedral angles summing up to $\pi$ to triples of
edges asymptotic to toric components of $\partial M$. If the
consistency and completeness equations are satisfied (see below),
this leads to a hyperbolic structure on $M$ minus $\alpha$ and the
toric components of $\partial M$.

\emph{Equations}. If an internal edge with non-zero dihedral angle
ends in a cusp then its length is infinity, so some of the length
equations must be dismissed when there are cusps. There are no
consistency issues connected with half-infinite edges but, when an
edge is infinite at both ends, one must make sure that the gluings
around the edge do not induce a sliding along the edge, which
translates into the condition that the \emph{similarity
moduli}~\cite{BePe:book} of the Euclidean triangles around the
edge have product $1$. This ensures existence of the hyperbolic
structure, but one still has to impose completeness of cusps. Just
as in the case where there are cusps only, this amounts to
requiring that the similarity tori on the boundary be Euclidean,
which translates into the \emph{holonomy equations} involving the
similarity moduli. Note that there is no completeness issue
connected to annular cusps.

\emph{Canonical decomposition}. When there are cusps, the set of
points $\calP$ to take the convex hull of consists of the norm-$1$
duals of the boundary components of the universal cover and of
some points on the light-cone dual to Margulis neighbourhoods of
the cusps. The precise discussion of how to choose these extra
points is too complicated to be reproduced here
(see~\cite{FriPe}), but in practice one has that the choice of
\emph{sufficiently small} Margulis neighbourhoods always leads to
the right result.

\section{Volume computation}\label{volume:section}
All the hyperbolic manifolds found by computer in~\cite{FriMaPe3}
can be decomposed into genuine tetrahedra, \emph{i.e.} tetrahedra
with positive dihedral angles (with ultra-ideal vertices, and
possibly some ideal ones). We do not know if all hyperbolic
manifolds admit such a genuine decomposition, but we know that if
we subdivide the Kojima decomposition into tetrahedra we find some
genuine tetrahedra and possibly some flat ones, which do not
contribute to the volume. Therefore the problem of computing the
volume of a hyperbolic manifold is completely reduced to the same
problem for a single tetrahedron with arbitrarily assigned
dihedral angles. We discuss in this section various formulae for
this volume. In doing this we will employ a unified viewpoint,
which includes finite, ideal, and ultra-ideal vertices.

\paragraph{Volumes of polyhedra: historical remarks}
The calculation of the volume of a polyhedron in 3-dimensional
space is a very old and difficult problem. The first known result
in this direction belongs to Tartaglia (1494) who found a formula
for the volume of a Euclidean tetrahedron, now known as the
Cayley-Menger formula. It was recently shown in~\cite{Sabitov}
and~\cite{CSW} that the volume of any Euclidean polyhedron is a
root of an algebraic equation whose coefficients are polynomial
functions of the lengths of the edges, the polynomials being
determined by the combinatorial type of the polyhedron.

In the hyperbolic and spherical spaces the situation is much more
complicated. We confine here to the hyperbolic case, and we first
recall that the volume formula for a biorthogonal tetrahedron
(\emph{orthoscheme}) has been known since Lobachevsky and
Schl\"{a}fli (see~\cite{Lb} and~\cite{Sh} respectively). The
volumes of the Lambert cube and of some other polyhedra were
computed by Kellerhals~\cite{Keller}, Derevnin and Mednykh
\cite{DM:Lambert}, Mednykh, Parker, and Vesnin \cite{MPV},
 and others. The volume of a regular tetrahedron in hyperbolic
space was investigated by Martin~\cite{Martin}. The volume formula for a
hyperbolic tetrahedron with a few non-ideal vertices was found by
Vinberg~\cite{Vi}.

Despite these partial results, a formula for the volume of an
arbitrary hyperbolic tetrahedron has been unknown until very
recently. The general algorithm for obtaining such a formula was
indicated by Hsiang ~\cite{H}, and the complete solution of the
problem was given within the space of a few years by several
authors, namely Cho and Kim ~\cite{ChK}, Murakami and
Yano~\cite{MY}, and Ushijima ~\cite{Ushijima:vol}.  In all these
papers the volume of a tetrahedron is expressed as an analytic
formula involving 16 dilogarithm or Lobachevsky functions whose
arguments depend on the dihedral angles of the tetrahedron and on
some additional parameter which is found as a root of some
complicated quadratic equation with complex coefficients.

A geometric meaning of the Murakami-Yano formula was recognized by
Leibon \cite{L} from the viewpoint of the so-called Regge
symmetry. An excellent exposition of these ideas and a complete
geometric proof of the Murakami-Yano formula was given by Mohanty
~\cite{M}. It is worth mentioning that the ideas of Regge symmetry
and scissors congruence were also partially used by Cho and Kim
in~\cite{ChK}, where the first general formula was actually
obtained.

A remarkable phenomenon is that the volume formulae for the
tetrahedron become much easier for \emph{symmetric} tetrahedra,
\emph{i.e.} for tetrahedra with identical dihedral angles at
opposite edges. This fact was first noticed by
Milnor~\cite{Milnor}, who expressed the volume of an ideal
tetrahedron (which is automatically symmetric) as the sum of the
values of the Lobachevsky function on the three dihedral angles.
It was later shown by Derevnin, Mednykh, and Pashkevich~\cite{DMP}
that a rather simple formula also exists for an arbitrary
symmetric tetrahedron with finite vertices.

In the next paragraph we will present an elementary integral
formula for the volume of an arbitrary tetrahedron in hyperbolic
space. The formula involves some parameters depending on the
dihedral angles, and it is very helpful to actually evaluate the
volume by computer.  The Murakami-Yano result can be obtained as
an easy consequence of this formula.

\paragraph{Volume of a tetrahedron without truncation}
Let $T(A,B,C,D,E,F)$ denote the tetrahedron in hyperbolic
$3$-space with dihedral
angles $A$, $B$, $C$, $D$, $E$, $F$ as in Fig.~\ref{tetra:fig}.
    \begin{figure}
    \begin{center}
    \input{tetra.pstex_t}
    \mylegend{Notation for the dihedral angles of a tetrahedron.} \label{tetra:fig}
    \end{center}
    \end{figure}
Recall~\cite{FriPe} that $T$ is determined up to isometry by $A$,
$B$, $C$, $D$, $E$, $F$. Moreover if the sum of the triple of
angles at a certain vertex is $\pi$ then the vertex is ideal.
Similarly, if the sum is less than $\pi$ then the vertex is
ultra-ideal, \emph{i.e.} the tetrahedron is truncated at the
vertex, which is then replaced by a triangle.

The next result is due to Derevnin and Mednykh~\cite{DM}:

\begin{teo}\label{volume:teo}
Suppose that all the vertices of $T=T(A,B,C,D,E,F)$ are ideal or
finite. Then ${\rm Vol}(T)$ is given by
$$-\frac14 \int\limits_{z_1}^{z_2} \log\frac
{\cos{\frac{A+B+C+z}{2}} \cdot \cos{\frac{A+E+F+z}{2}}\cdot
\cos{\frac{B+D+F+z}{2}}\cdot \cos{\frac{C+D+E+z}{2}}}
{\sin{\frac{A+B+D+E+z}{2}}\cdot \sin{\frac{A+C+D+F+z}{2}}\cdot
\sin{\frac{B+C+E+F+z}{2}}\cdot \sin{\frac{ z}{2}}} \,{\rm d}z , $$
with $z_1$ and $z_2$ given by
$$z_1=\arctan\frac{{k_2}}{{k_1}}-\arctan\frac{k_4}{k_3},\qquad
z_2=\arctan\frac{{k_2}}{{k_1}}+\arctan\frac{k_4}{k_3},$$
where
\begin{eqnarray*}
{k_{1}} &=&-(\cos S+\cos (A+D)+\cos (B+E)+\cos (C+F)\\
&&+\cos(D+E+F) +\cos(D+B+C)\\
&&+\cos (A+E+C)+\cos (A+B+F)), \\
{k_{2}} &=&\sin S+\sin (A+D)+\sin (B+E)+\sin (C+F)\\
&&+\sin (D+E+F)+\sin (D+B+C)\\
&&+\sin (A+E+C)+\sin (A+B+F), \\
{k}_{3} &=&2(\sin {A}\sin {D}+\sin {B}\sin {E}+\sin {C}\sin{F}),\\
{k}_{4} &=&\sqrt{k_1^2+k_2^2-k_3^2},
\end{eqnarray*}
and $S=A+B+C+D+E+F$. Moreover the $k_j$'s and $z_j$'s are all real
numbers, so the integral is just an ordinary integral on an
interval of the real line, and the function to be integrated
vanishes at the $z_j$'s.
\end{teo}

\begin{rem}\label{vertices:rem}
\emph{There is a very transparent geometric interpretation of the
sums of dihedral angles $V_1=A+B+C$, $V_2=A+E+F$, $V_3=B+D+F$, and
$V_4=C+D+E$ appearing in the numerator of the volume formula, and
of the sums $H_1=A+B+D+E$, $H_2=A+C+D+F$, and $H_3=B+C+E+F$
appearing in the denominator. Namely, the $V_j$'s correspond to
the triples of edges incident to the vertices of the tetrahedron,
while the $H_j$'s correspond to the Hamiltonian cycles.}
\end{rem}

\begin{rem}\label{equation:rem}
\emph{The parameters $z_1$ and $z_2$ appearing in the volume
formula can be shown to be roots of the equation ${k_1}\,\cos z
+{k_2}\,\sin z={k_3},$ while ${k}_{4}^2=-4\det(G),$ where $G$ is
the Gram matrix of $T$. The numbers $z_1$ and $z_2$ also have a
geometric meaning, as explained in~\cite{M}. Namely, they arise as
parameters for the decomposition of an ideal octahedron into four
ideal tetrahedra with an edge in common. The octahedron is
canonically defined by the tetrahedron $T$, and its dihedral
angles are just linear combinations of those of $T$.}
\end{rem}

Recall now that the dilogarithm function is defined by the
integral
$${\rm Li}_2(x)=-\int\limits_0^x\frac{\log(1-t)}{t}\,{\rm d}t,$$ where $x\in\matC \setminus{[}1,\infty{)}$
and $\log$ is the continuous branch of the logarithm function
given by $\log \xi = \log|\xi|+i\,\arg \xi$ with the constraint
$-\pi< \arg \xi< \pi$. Let
$$l(z)= {\rm Li}_2\,({\rm e}^{i\,z}).$$
As an immediate consequence of Theorem~\ref{volume:teo} we have
the following Murakami-Yano-Ushijima formula obtained in~\cite{MY}
and~\cite{Ushijima:vol}:

\begin{cor}\label{dilogarithm:cor}
Suppose that all the vertices of $T=T(A,B,C,D,E,F)$ are ideal or
finite. Then
$${\rm Vol}(T)=\frac{1}{2}\Im\Big(U(z_1,\,T)-U(z_2,\,T)\Big),$$
where
\begin{eqnarray*}
U(z,T)&=&\frac12\Big(l(z)+l(A+B+D+E+z)\\
&&\quad +l(A+C+D+F+z)+l(B+C+E+F+z)\\
&&\quad -l(\pi+A+B+C+z)-l(\pi+A+E+F+z)\\
&&\quad -l(\pi+B+D+F+z)-l(\pi+C+D+E+z)\Big). \\
\end{eqnarray*}
\end{cor}

Conerning this formula we note that
$$\Im( l(z))=\Im( {\rm Li}_2({\rm e}^{i\,z}))=2\Lambda(z/2),$$ 
where
$\Lambda(z)$ is the Lobachevsky function defined by the integral
$$\Lambda(\theta)=-\int\limits_0^\theta \log|2\sin t|\,{\rm d}t.$$
This shows, as already announced, that the volume of a
tetrahedron is an algebraic sum of 16 Lobachevsky functions.

\paragraph{Truncated tetrahedra}
According to~\cite{Ushijima:vol} and~\cite{Vi}, the determinant of the
Gram matrix is strictly negative also for truncated tetrahedra, so the numbers 
$k_1,\ldots,k_4$ appearing in the statement of Theorem~\ref{volume:teo} are still real.
However a non-trivial issue arises concerning the choice of the analytic branch
of the $\arctan$ function for the definition of the numbers $z_1$ and $z_2$.
More precisely, one can be forced to change branch, in order to ensure the continuity of $z_1$ and $z_2$,
when $k_1$ approaches $0$, which indeed can happen when there
are some ultra-ideal vertices. For example, for
$$(A,B,C,D,E,F)=\left(\frac{\pi}{12},\frac{\pi}{3},\frac{\pi}{10.18},\frac{\pi}{12},\frac{\pi}{3},\frac{\pi}{10.18}\right)$$
we have $k_1 > 0$, whereas for
$$(A,B,C,D,E,F)=\left(\frac{\pi}{12},\frac{\pi}{3},\frac{\pi}{10.19},\frac{\pi}{12},\frac{\pi}{3},\frac{\pi}{10.19}\right)$$
we obtain  $k_1 < 0$.  Therefore in this case we have to change the analytic branch of
$\arctan (k_2/k_1)$. The most convenient way to do so is to replace $\arctan(k_2/k_1)$ by $\pi/2-\arctan(k_1/k_2)$,
which yields a continuous variation of the values of $z_1$ and $z_2$ as $k_1$ passes through 0. We note that $k_3$ and $k_4$
are always positive, so there is no such problem with $\arctan(k_4/k_3)$. One can now check that, after 
appropriately choosing the analytic branches of the $\arctan$ function, the formula of 
Theorem~\ref{volume:teo} still gives the correct value of the volume for truncated tetrahedra.

\paragraph{Symmetric tetrahedra}
We will say that $T=T(A,B,C,D,E,F)$ is \emph{symmetric} if $A=D$,
$B=E$, and $C=F$, and in this case we will denote $T$ just by
$T(A,B,C)$ for simplicity.

We begin with the easy case of ideal tetrahedra, which are
automatically symmetric. Moreover $T(A,B,C)$ is ideal if and only
if $A+B+C=\pi$. In this case the volume turns out~\cite{Milnor} to
be given simply by
$${\rm Vol}(T(A,B,C))=\Lambda(A)+\Lambda(B)+\Lambda(C),$$ where $\Lambda$ is the
Lobachevsky function defined above.
For symmetric tetrahedra $T(A,B,C)$ with finite vertices,
\emph{i.e.} such that $A+B+C>\pi$, the following was shown
in~\cite{DMP}:

\begin{teo}\label{symvolumehyp}
If $T=T(A,B,C)$ has finite vertices then ${\rm Vol}(T)$ is given by
$$2\int\limits_{\theta}^{\pi/2} \frac
{\sin^{-1} (\cos A \cos t)+\sin^{-1} (\cos B \cos t)
+\sin^{-1} (\cos C \cos t)-\sin^{-1} (\cos t)}{\sin 2t}{\rm d}t$$
where $\theta\in(0,\pi/2)$ satisfies
$$
\tan \theta= \frac{1-a^{2}-b^{2}-c^{2}-2abc}
{\sqrt{(1-a+b+c)(1+a-b+c)(1+a+b-c)(-1+a+b+c)}},
$$
with $a=\cos A$, $b=\cos B$, and $c=\cos C$.
\end{teo}

\begin{rem}\label{sinerule:rem}
\emph{The value of $\theta$ in the previous statement has a simple
geometric interpretation in terms of the ``sine rule''
(see~\cite[Theorem 7]{DMP}): if $l_{A},l_{B},l_{C}$ are the
lengths of the edges of $T$ with dihedral angles $A,B,C$
respectively, then}
$$\frac{\sin A}{\sinh l_{A}}\,=\,\frac{\sin B}{\sinh l_{B}}
\,=\,\frac{\sin C}{\sinh l_{C}}\,=\tan\theta.$$
\end{rem}

We conclude this section by noting that no simple formula is
currently known for the volume of a symmetric tetrahedron with
ultra-ideal vertices, but it is reasonable to expect that such a
formula should exists and have connections with an ultra-ideal
version of the sine rule.

\section{Manifolds without annular cusps}\label{cpt:section}
In this section we recall the definition and main properties of a
certain class of manifolds $\calM_{g,k}\subset\calH_{g+k}$ studied
in~\cite{FriMaPe1, FriMaPe2}, and we describe the experimental
results of~\cite{FriMaPe3} on $\calH_n$ for $n\leqslant 4$. Recall
that $\calH_n$ contains the hyperbolic manifolds of complexity $n$
with non-empty \emph{compact} geodesic boundary. The next section
will be devoted to manifolds with non-compact boundary,
\emph{i.e.} with annular cusps.

\paragraph{A special class of manifolds}
Let us denote by $\Sigma_g$ the closed orientable surface of genus
$g$. For $g\geqslant 2$ and $k\geqslant 0$ we define $\calM_{g,k}$
as the set of all compact orientable manifolds $M$ having an ideal
triangulation with $g+k$ tetrahedra, and $$\partial M =
\Sigma_g\sqcup\left(\mathop{\sqcup}\limits_{i=1}^k
T_i\right)\qquad {\rm with}\ T_i\cong\Sigma_1.$$ As a motivation
for this definition, we mention here that an ideal triangulation
of a manifold whose boundary is the union of $\Sigma_g$ and $k$
tori contains at least $g+k$ tetrahedra. So $\calM_{g,k}$ is the
set of manifolds $M$ having the smallest possible complexity,
given the topological constraints on $\partial M$.

Our first result shows that the class just introduced is very
large:

\begin{prop}
\begin{itemize}

\item $\calM_{g,k}$ is non-empty precisely for $g>k$ or $g=k$ and
$g$ even;

\item The values of $\#\calM_{g,k}$ for small $g$ and $k$ are
as shown in Table 1;

\begin{table}
\begin{center}
\makebox{
\begin{tabular}{c||c|c|c}
\phantom{\Big|}                     & $k=0$     & $k=1$     & $k=2$     \\
\hline\hline \phantom{\Big|} $g=2$  & $8$       & $1$       & $1$       \\
\hline \phantom{\Big|} $g=3$        & $74$      & $12$      & $1$       \\
\hline \phantom{\Big|} $g=4$        & $2340$    & $416$     & $51$      \\
\hline \phantom{\Big|} $g=5$        & $97568$   & $17900$   & ?         \\
\end{tabular}
}
\end{center}
\mylegend{Some values of $\#\calM_{g,k}$.}
\end{table}

\item For any fixed $k$ there exist constants $C>c>0$ such that
$$g^{c\cdot g} < \#\calM_{g,k} < g^{C\cdot g}\qquad {\it for}\qquad g\gg 0.$$

\end{itemize}
\end{prop}

We now give the main statement we have about the elements of
$\calM_{g,k}$ and their topological and geometric invariants. We
address the reader to~\cite{scharlemann} for the definition of the
\emph{Heegaard genus} of  a triple $(M,\partial_0M,\partial_1M)$,
and to~\cite{TuVi} for the definition of the Turaev-Viro
invariants $\TV_r(M)$ of a compact $3$-manifold $M$, for
$r\geqslant 2$.

\begin{teo} \label{main:teo}
Let $M\in\calM_{g,k}$. The following holds:

\begin{enumerate}

\item \label{main:hyperbolic:item} $M\setminus(\sqcup T_i)$ is
hyperbolic, and its volume depends only on $g$ and $k$;

\item \label{main:canonical:item} $M$ has a unique ideal
triangulation with $g+k$ tetrahedra, which gives the canonical
Kojima decomposition of $M\setminus(\sqcup T_i)$;

\item \label{main:complexity:item} $M$ has complexity $g+k$;

\item \label{main:genus:item} The Heegaard genus of
$\left(M,\Sigma_g,\sqcup T_i\right)$ is $g+1$;

\item \label{main:homology:item} $H_1(M;\matZ)=\matZ^{g+k}$;

\item \label{main:turaev:viro:item} The Turaev-Viro invariant
$\TV_r(M)$ of $M$ depends only on $r$, $g$ and $k$.

\end{enumerate}
\end{teo}

The first two points of this theorem are established using
precisely the philosophy of moduli, equations, and tilts sketched
in Section~\ref{equations:section}. Namely, given a minimal
triangulation of $M\in\calM_{g,k}$, one shows that each
tetrahedron has either one or no vertex asymptotic to a toric
cusp. If there is one such vertex, one chooses the dihedral angles
to be $\pi/3$ at the edges ending at that vertex, and to be all
equal to some $\alpha$ at the other edges. If there is no such
vertex, one chooses the tetrahedron to have all dihedral angles
equal to some $\beta$. Using a continuity argument then one sees
that there exist values of $\alpha$ and $\beta$ satisfying the
consistency equations. Completeness and the computation of tilts
is then straight-forward.

We also notice that the previous results show the tremendous power
of hyperbolic geometry compared to the topological invariants: the
classes $\calM_{g,k}$ are extremely large, and only hyperbolic
geometry is able to distinguish their elements from each other.

We refrain from recalling the precise statements and details here,
but we want to mention that a thorough analysis of the Dehn
fillings of the elements of $\calM_{g,k}$ was carried out
in~\cite{FriMaPe3}, leading in particular to the solution of a
problem raised by Gordon and Wu~\cite{GoWu-AA, GoWu-DA, Wu:1, Wu:2}
on the maximal distance between a
boundary-reducible and a non-acylindrical slope on a ``large''
hyperbolic manifold.

\paragraph{Experimental results}
We describe here $\calH_n$ for $n\leqslant 4$. 
It is easy to see that $\calH_1$ is empty.
Moreover, $\calH_2$ was shown in~\cite{Fujii} to have $8$
elements, all with the same volume $\approx 6.451990$. Therefore
$\calH_2$ coincides with the set $\calM_{2,0}$ discussed above, so
all its members share the same invariants, except the hyperbolic
structure itself.

The strategy mentioned in Section~\ref{complexity:section} has
been implemented in~\cite{FriMaPe2} to classify $\calH_3$ and
$\calH_4$, leading to the following results. We denote here by
$\calK(Y)$ the Kojima canonical decomposition of a hyperbolic $Y$,
and we recall that $\calHtilde_n$ is the set of candidate
hyperbolic $3$-manifolds of complexity $n$. We emphasize that the
next results indeed have an experimental nature, but their
validity was confirmed by a number of computations by hand and
cross-checks. We also mention that all the values of the volume
are approximate ones. More accurate approximations are available
on the web~\cite{www}.

\begin{teo}\label{H3:teo}
\begin{itemize}

\item $\calH_3$ coincides with $\widetilde{\calH}_3$ and it has
$151$ elements;

\item $\calM_{3,0}\subset\calH_3$ consists of $74$ elements of
volume $10.428602$;

\item $\calM_{2,1}\subset\calH_3$ consists of a single manifold of
volume $7.797637$;

\item The $76$ elements of $\calH_3\setminus(\calM_{3,0}\cup
\calM_{2,1})$ all have boundary $\Sigma_2$, and they split as
follows:

\begin{itemize}

\item $73$ compact $Y$'s with $\calK(Y)$ consisting of three
tetrahedra; the volume function attains on them $15$ different
values ranging between $7.107592$ and $8.513926$, with maximal
multiplicity $9$;

\item $3$ compact $Y$'s with $\calK(Y)$ consisting of four
tetrahedra; they all have the same volume $7.758268$.

\end{itemize}
\end{itemize}
\end{teo}

\begin{teo}
\begin{itemize}

\item $\calH_4$ has $5033$ elements, and $\widetilde{\calH}_4$ and
has $6$ more;

\item $\calM_{4,0}\subset\calH_4$ has $2340$ elements of volume
$14.238170$;

\item $\calM_{3,1}\subset\calH_4$ has $12$ elements of volume
$11.812681$;

\item $\calM_{2,2}\subset\calH_4$ has a single element of volume
$9.134475$;

\item The $2680$ elements $Y$ of $\calH_4$ not belonging to any
$\calM_{g,4-g}$ split as described in Table 2 according to the
boundary of $\Ybarbar$ (columns) and type of blocks of the Kojima
decomposition (rows). Some volume information has also been
inserted in each box, namely the unique value of volume if there
is one, or the minimum and the maximum of volume, the number of
values it attains, and the maximal multiplicity of these values.

\end{itemize}

\end{teo}

\begin{table}
\begin{center}
\begin{tabular}{l||c|c|c}
& $\Sigma_3$ & $\Sigma_2$ &$\Sigma_2$,\ 1\ cusp\\ \hline\hline 4
tetrahedra &
    {\bf 1936} & {\bf 555} & {\bf 16} \\
    & $\min=11.113262$ & $\min=7.378628$ & $\min=8.446655$ \\
    & $\max=12.903981$ & $\max=10.292422$ & $\max=9.774939$ \\
    & ${\rm values}=59$ & ${\rm values}=169$ & ${\rm values}=8$ \\
    & ${\rm max\ mult}=138$ & ${\rm max\ mult}=27$ & ${\rm max\ mult}=3$ \\ \hline
5 tetrahedra & {\bf 42} & {\bf 41} & \\
    & ${\rm vol}=11.796442$ & $\min=8.511458$ & \\
    & & $\max=9.719900$ & \\
    & & ${\rm values}=16$ & \\
    & & ${\rm max\ mult}=6$ & \\ \hline
6 tetrahedra & & {\bf 3} & \\
    & & ${\rm vol}=8.297977$ & \\ \hline
8 tetrahedra & & {\bf 3} & \\
    & & ${\rm vol}=8.572927$  & \\ \hline
1 octahedron & {\bf 56} & {\bf 14} & \\
(regular) & ${\rm vol}=11.448776$ & ${\rm vol}=9.415842$ & \\
\hline
1 octahedron & & {\bf 8} & \\
(non-regular) & & ${\rm vol}=8.739252$ & \\ \hline
2 pyramids & & {\bf 4} & {\bf 2} \\
with square & & ${\rm vol}=9.044841$ & ${\rm vol}=8.681738$ \\
basis & & &
\end{tabular}\end{center}
\mylegend{Elements of $\calH_4$ not belonging to any $\calM_{g,4-g}$.}
\end{table}

We conclude this section with some information on the actual
computer implementation of the enumeration process. We solve the
hyperbolicity equations using Newton's method with partial
pivoting, after explicitly writing the derivatives of the length
function. Convergence to the solution is always extremely fast,
and it can be checked to be stable under modifications of the
numerical parameters involved in the implementation of Newton's
method. Concerning the Kojima decomposition, we mention that the
evolution of a triangulation toward the canonical decomposition is
not quite sure to converge in general, but it always does in
practice. We also point out that our computer program is only able
to handle triangulations: whenever some mixed negative and zero
tilts appear, the canonical decomposition must be worked out by
hand and actually proved not to be a triangulation.

\section{Manifolds with annular cusps}\label{CFMP:section}
Given a hyperbolic manifold $Y$ with annular cusps, one may decide
to use either $c(\Ybar)$ of $c(\Ybarbar)$ as a definition of the
complexity of $Y$, the former being more natural from a
topological viewpoint, and the latter from a geometric viewpoint.
However, as announced in Proposition~\ref{bad:annuli:prop} and
explained in this section, neither definition allows to employ the
powerful techniques of special spines and ideal triangulations to
actually carry out the enumeration of manifolds in order of
increasing complexity. For this reason the understanding of
manifolds with annular cusps is still very limited. This section
is only devoted to the description of a special class of such
manifolds, but the properties of this class are already sufficient
to show that the set of hyperbolic manifolds with annular cusps is
very large, and that there is very little control on the topology
of the compactifications of such manifolds. We will also see that
the geometric information on this special class allows to prove
some very interesting (and apparently unrelated) results. We
address the reader to~\cite{CFMP} for all proofs and further
details.

\paragraph{Combinatorial triangulations}
We denote by $\calT_n$ the set of all possible simplicial pairings
between the faces of $n$ tetrahedra.  We view two pairings to be
the same if they coincide combinatorially, and we note that in
general a pairing can fail give an ideal triangulation of a
manifold, because the link of the midpoint of an edge can be the
projective plane. Moreover, if a pairing actually gives a
manifold, then this manifold need not be orientable. However, if
we fix an orientation on the tetrahedron and require all the
pairing maps to reverse the orientation, the result is an ideal
triangulation of an orientable manifold. We denote by $\calT$ the
union of all $\calT_n$'s.

Given $T\in\calT$ we define a space $Y(T)$ by performing the
pairings in $T$ and then removing first an open regular
neighbourhood of the vertices and then a closed regular
neighbourhood of the edges.

\begin{rem}\label{YTbar:rem}
\emph{The manifold $Y(T)$ is a non-compact one with boundary, and
its topological ends have the shape of the product of an annulus
and a closed half-line. Let us define $\Ybar(T)$ and $\Ybarbar(T)$
as the natural compactifications of $Y(T)$, as we did in
Section~\ref{structure:section} for a hyperbolic $Y$. Then
$\Ybar(T)$ is a (possibly non-orientable) handlebody, with genus
$n+1$ if $T\in\calT_n$. Moreover, if $T$ is an ideal triangulation
of a manifold $M$, then $\Ybarbar(T)$ is homeomorphic to $M$.}
\end{rem}

\paragraph{Relative handlebodies and hyperbolicity}
For $T\in\calT$ the manifold $Y(T)$ can be viewed as
$H\setminus\Gamma$, where $H$ is a handlebody and
$\Gamma\subset\partial H$ is a system of disjoint loops. At the
risk of a little ambiguity, which eventually will not create any
serious problem, we then identify $Y(T)$ to the pair $(H,\Gamma)$,
and, following Johannson~\cite{johannson}, we call such a pair a
\emph{relative handlebody}. We denote by $\calA_n$ the set of all
$Y(T)$'s as $T$ varies in $\calT_n$, and by $\calA$ the union of
all $\calA_n$'s.

To state our first result we call \emph{complexity} of a relative
handlebody $(H,\Gamma)$ the minimum of $|\Gamma\cap\partial D|$
where $D$ is a system of disjoint properly embedded discs in $H$
such that $\partial D$ cuts $\partial H$ into a union of pairs of
pants.

\begin{prop}\label{Tn:prop}
\begin{itemize}
\item For all $T\in\calT$ the relative handlebody $Y(T)$ is
hyperbolic. \item The map $T\mapsto Y(T)$ gives a bijection
between $\calT$ and $\calA$. \item Among the hyperbolic
$(H,\Gamma)$'s of genus $n+1$, the elements of $\calA_n$ can be
characterized as those having minimal complexity, equal to
$10\cdot n$, and as those having minimal volume, equal to $n\cdot
v_O$, where $v_O\approx 3.66386$ is the volume of a hyperbolic
regular ideal octahedron.
\end{itemize}
\end{prop}

The proof of this result entirely depends on the techniques
described in Section~\ref{equations:section}. One first notices
that a regular ultra-ideal tetrahedron with all dihedral angles
equal to $0$ is actually a regular ideal octahedron with a
checkerboard coloring of the faces, the white ones being internal
faces and the black ones being truncation faces. Then one sees
that for each given $T\in\calT_n$ the corresponding $Y(T)$ can be
constructed by gluing according to $T$ the white faces of $n$
regular ideal octahedra. To conclude, one proves that $T$ itself
is the canonical Kojima decomposition of $Y(T)$.

Remark~\ref{YTbar:rem} and Proposition~\ref{Tn:prop} are already
sufficient to prove Proposition~\ref{bad:annuli:prop} which shows,
as already mentioned, that there are very many manifolds with
annular cusps and that their topology is rather arbitrary. The
next result goes in the same direction. We call \emph{tangle} in a
compact 3-manifold with boundary a finite union of disjoint
properly embedded arcs.

\begin{cor}\label{tangles:enlarged:intro:cor}
Every tangle in every compact $3$-manifold is contained in a
tangle whose complement lies in $\calA_n$ for some $n$ (in
particular, it is hyperbolic with annular cusps).
\end{cor}

An unexpected consequence of Proposition~\ref{Tn:prop} is also the
following:

\begin{teo}\label{triangulations:determined:intro:prop}
Let $T_0$ and $T_1$ be triangulations of the same compact
$3$-manifold, with finite and/or ideal vertices, possibly with multiple and
self-adjacencies. Assume that the $1$-skeleta of $T_0$ and $T_1$
coincide. Then $T_0$ and $T_1$ are isotopic relatively to the
1-skeleton.
\end{teo}

\paragraph{Doubles and Dehn filling}
We state in this paragraph some of the surprising results one can
establish starting from the construction of the hyperbolic
relative handlebodies $Y(T)$.

For $T\in\calT$, we define $D(T)$ to be the ``orientable double'' of $Y(T)$,
namely either the union of two copies of $Y(T)$ along the boundary, when $Y(T)$
is orientable, or the quotient of the orientation covering of $Y(T)$ under the
restriction to the boundary of the involution, when $Y(T)$ is non-orientable.
We define $\calD$ as the set of all $Y(T)$ for $T\in\calT$. The next result
shows that the family $\calD$ is completely
classified in combinatorial terms and it is universal for $3$-manifolds
under the operation of Dehn filling:

\begin{teo}\label{T:to:D:intro:teo}
\begin{itemize}
\item Every member of $\calD$ is an orientable cusped hyperbolic
manifold without boundary, and it is the complement of a link in a
connected sum of some copies of $S^2\times S^1$; \item The
correspondence $T\mapsto D(T)$ defines a bijection between $\calT$
and $\calD$. \item Every closed orientable $3$-manifold is a Dehn
filling of a manifold in $\calD$.
\end{itemize}
\end{teo}

Since a Dehn filling of a hyperbolic manifold ``typically'' is
hyperbolic, the class $\calD$ provides a powerful method to
construct closed hyperbolic manifolds. As an application of this
method one can establish the following refinement of the main
result of~\cite{Kojima:isom}:

\begin{prop}\label{all:groups:intro:teo}
There exists $c>0$ such that, given a
finite group $G$, there is a closed orientable
hyperbolic $3$-manifold $M$ with ${\rm Isom}(M)\cong G$ and ${\rm Vol}(M)\leqslant c\cdot |G|^9$.
\end{prop}

We now recall that,
according to Thurston's hyperbolic Dehn filling theorem, on each
cusp of a finite-volume hyperbolic 3-manifold there is only a
finite number of slopes filling along which one gets a
non-hyperbolic 3-manifold. These slopes are called
\emph{exceptional}, and a considerable effort has been devoted to
understanding them~\cite{Gordon}. If $T$ is a triangulation, the
hyperbolic manifold $D(T)$ has a preferred horospherical cusp
section, and each component of this section corresponds to an edge
of $T$. Moreover the valence of the edge gives a lower bound for
the length of the second shortest geodesic on the component. This
fact and the Agol-Lackenby 6-theorem~\cite{Agol:length-6, Lackenby}
imply the following:

\begin{prop}\label{exceptional:intro:prop}
If every edge of $T$ has valence at least $7$ then there is at
most one exceptional slope on each cusp of $D(T)$.
\end{prop}

We conclude this paper by showing that in some cases the
combinatorics of a triangulation is already sufficient to prove
hyperbolicity of the underlying manifolds. To motivate our result,
suppose first that a 3-manifold $M$ has an \emph{ideal}
triangulation $T$ in which each edge has valence at least $6$. An
easy argument shows that $\chi(T)\leqslant 0$, and that
$\chi(T)=0$ precisely when all valences are $6$. Moreover, in the
last case, the boundary of $M$ is a disjoint union of tori and
Klein bottles, and Thurston's hyperbolicity equations for cusped
manifolds have a very simple solution, given by regular ideal
tetrahedra. Analogously, if all edges have one and the same
valence $v\geqslant 7$, the hyperbolicity equations for the
geodesic boundary case have a simple solution, given by regular
truncated tetrahedra with dihedral angles $2\pi/v$. An argument
based on the Agol-Lackenby machinery~\cite{Agol:length-6,
Lackenby} allows to generalize these facts as follows:

\begin{prop}
If $M$ has an ideal triangulation $T$ whose edges have valence at
least $6$, then $M$ is hyperbolic, and the edges of $T$ are
homotopically non-trivial relative to $\partial M$.
\end{prop}

According to the last assertion of this result, the edges of $T$
can be straightened to geodesics in $M$ without entirely disappearing
into infinity, which suggests that $T$ itself can be straightened.
We address the reader to~\cite{ideal:surv} for more information
on the general problem of existence of a triangulation which can be
straightened to an ideal one.

\vspace{1.5cm}

\noindent
Sobolev Institute of Mathematics\\
Russian Academy of Sciences\\
4 Acad. Koptyug avenue, Novosibirsk 630090, Russia\\
smedn@mail.ru

\vspace{.5cm}

\noindent
Dipartimento di Matematica Applicata\\
Universit\`a di Pisa\\
Via Bonanno Pisano 25B, 56126 Pisa, Italy\\
petronio@dm.unipi.it

\end{document}

%% file: tetra.pstex_t
\begin{picture}(0,0)%
\includegraphics{tetra.pstex}%
\end{picture}%
\setlength{\unitlength}{3355sp}%
\begingroup\makeatletter\ifx\SetFigFont\undefined%
\gdef\SetFigFont#1#2#3#4#5{%
  \reset@font\fontsize{#1}{#2pt}%
  \fontfamily{#3}\fontseries{#4}\fontshape{#5}%
  \selectfont}%
\fi\endgroup%
\begin{picture}(2949,3324)(3364,-5173)
\put(3833,-2918){\makebox(0,0)[lb]{\smash{\SetFigFont{10}{12.0}{\familydefault}{\mddefault}{\updefault}{$A$}%
}}}
\put(3773,-4381){\makebox(0,0)[lb]{\smash{\SetFigFont{10}{12.0}{\familydefault}{\mddefault}{\updefault}{$B$}%
}}}
\put(4231,-3428){\makebox(0,0)[lb]{\smash{\SetFigFont{10}{12.0}{\familydefault}{\mddefault}{\updefault}{$C$}%
}}}
\put(5648,-4358){\makebox(0,0)[lb]{\smash{\SetFigFont{10}{12.0}{\familydefault}{\mddefault}{\updefault}{$D$}%
}}}
\put(5708,-2528){\makebox(0,0)[lb]{\smash{\SetFigFont{10}{12.0}{\familydefault}{\mddefault}{\updefault}{$E$}%
}}}
\put(4966,-4013){\makebox(0,0)[lb]{\smash{\SetFigFont{10}{12.0}{\familydefault}{\mddefault}{\updefault}{$F$}%
}}}
\end{picture}